\documentclass{amsart}
\usepackage{amssymb,latexsym}
\usepackage{epsfig}
\theoremstyle{plain}

\newtheorem{proposition}{Proposition}
\newcommand{\rn}{\mathbb{R}}
\newcommand{\join}{\cup}
\newcommand{\bigjoin}{\bigcup}
\newcommand{\meet}{\cap}

\renewcommand{\le}{\leqslant}
\renewcommand{\ge}{\geqslant}

\newcommand{\tchk}{\cdot}
\newcommand{\A}{\Tilde{A}}

\newcommand{\spann}{\text{span}}
\newcommand{\eig}{\text{eig}}
\newcommand{\per}{\text{per}}
\newcommand{\bez}{\backslash}
\newcommand{\supp}{\text{supp}}
\newcommand{\inter}{\text{int}}
\begin{document}
\title{Max-plus definite matrix closures and their eigenspaces}
\author{Serge\u{\i} Sergeev}
\address{Serge\u{\i} Sergeev, Sub-Department of Quantum Statistics and Field Theory, Department of
Physics, Moscow State University, 199992 Vorobievy Gory, Moscow, Russia}
\email{sergiej@gmail.com}
\thanks{Research supported by RFBR grant number 05-01-00824}
\begin{abstract}
In this paper we introduce the definite closure operation
for matrices with 
%%%AS 
finite
permanent, reveal inner structures
of definite eigenspaces, and establish some facts about Hilbert distances
between these inner structures and 
%%%AS
the
boundary of the definite eigenspace.
\end{abstract}
\keywords{Max-algebra, max-plus geometry, definite matrix, 
max-plus eigenspace, Hilbert distance}
\subjclass[2000]{15A18, 06F15, 52A99.}
\maketitle

\section{Introduction}\label{s:intro}

This paper is a contribution to geometrical understanding
of some algebraic results on max-plus eigenspaces that were obtained
by P.~Butkovi\v{c}
in~\cite{But00} (see also~\cite{But03}). 
The sources of geometrical inspiration
for our work are~\cite{CGQ04} and~\cite{LMS01}, as well as~\cite{DSS} 
and~\cite{DevSt}. Our approach is closer to 
that of~\cite{DSS} and~\cite{DevSt}, since we think of max-plus semiring
with its simplifying total order rather than of 
generalized algebraic structures of idempotent analysis. 
However, our viewpoint differs from that of~\cite{DSS} and~\cite{DevSt}
in that we use
%%%AS  
basic tools of max-algebra instead of the more
sophisticated machinery of convex geometry.

The paper is organized as follows. In Sect.~\ref{s:maxplus} we recall
the basic tools of max-algebra that we need. In Sect.~\ref{s:defclos}
we introduce 
%%%AS 
definite forms for max-plus matrices with nonzero (finite)
permanent and prove that closures of all definite
forms of a given matrix coincide. Thus we can introduce the `definite closure'
of any max-plus matrix with nonzero permanent. 
We also introduce 
%%%AS
definite eigenspaces, make 
%%%AS 
some observations 
on systems
of inequalities that define them, and consider an application
to the cellular decomposition introduced in~\cite{DevSt}. 
In Sect.~\ref{s:hdist} we use a
representation, due to~\cite{But00}, of the definite eigenspace that
reveals some inner structures. Then
we establish some facts about Hilbert distances between these inner structures
and 
%%%AS
the
boundary of the eigenspace.

The author wishes to thank P.~Butkovi\v{c}, A.~Churkin, A.~Kurnosov,  
A.~Sobolevski\u{\i}, and anonymous referees for helpful comments on the paper.

\section{Some tools of max-algebra}\label{s:maxplus}

In its ordinary setting, max-algebra is linear algebra
over the semiring $\rn_{\max}$. This semiring is the set 
$\rn\join\{-\infty\}$ equipped with the operations of `addition' 
$\oplus=\max$ and `multiplication' $\odot=+$. Its `zero' $\textbf 0$ is
equal to $-\infty$ and its `unity' $\textbf 1$ is equal to~$0$.

The semiring $\rn_{\max}$
resembles $\rn^+$, the positive part of the field of real numbers:
in 
%%%AS 
both structures multiplication admits inverses and enjoys
distributivity over addition, but subtraction is not allowed. 
In $\rn_{\max}$, however, the addition $\oplus$ is \emph{idempotent}, i.e.\
for any $\alpha\in\rn_{\max}$ we have $\alpha\oplus\alpha=\alpha$.
In $\rn^+$ this is certainly not the case.
The semiring $\rn_{\max}$ can be seen as an idempotent 
`dequantization' of $\rn^+$ (see, e.g.,
\cite{LitMas95} and \cite{KolMas97}). 

In max-algebra, it is also possible to 
%%%AS
exponentiate  
and to take roots.
These operations are nothing but conventional multiplication
and conventional division, respectively. Indeed, for any $\alpha\ne\textbf 0$
one has $\alpha^n=\alpha\times n$ and $\alpha^{\frac{1}{n}}=\alpha/n$ 
(for the remaining case $\textbf{0}^n=\textbf{0}^{\frac{1}{n}}=\textbf 0$ 
since we assume that $\alpha\odot\textbf{0}=\textbf 0$ for any $\alpha$)

%%%AS
One of the principal objects 
of max-algebra is $\rn_{\max}^n$, 
the set of $n$-component vectors with components
in $\rn_{\max}$. This set is equipped with `addition' 
$(x\oplus y)_i=x_i\oplus y_i$ and with `multiplication' by any 
\emph{max-plus scalar} (i.e. by any element of $\rn_{\max}$)
$(\alpha\odot y)_i=\alpha\odot y_i$. 
The set~$\rn_{\max}^n$ equipped with these operations, as well
as subsets of~$\rn_{\max}^n$ closed under these operations, will
be called \emph{max-plus spaces}.
Below the notation $\odot$ will be frequently omitted.

These max-plus spaces resemble
linear spaces in that the laws of associativity and distributivity
hold, but again there is no subtraction and there is idempotency of
addition. 
%%%AS
Structures  
of this kind are called \emph{idempotent
semimodules} and are 
%%%AS
a central object of the  
study in the \emph{idempotent
analysis}, see \cite{KolMas97}. 

A max-plus space 
%%%AS
$S\subset \rn_{\max}^n$  
is said to be \emph{finitely generated}
if there is a set of vectors $\{\,v_1,\ldots,v_s\,\}$ such that for any
$y\in S$ one can find scalars $\alpha_1,\ldots,\alpha_s$ such that
$y=\bigoplus_{i=1}^s\alpha_iv_i$. The set $\{\,v_1,\ldots,v_s\,\}$ is the
\emph{generating set} of $S$. It is \emph{minimal} if no $v_i$
can be expressed as a linear combination of the other generators, i.e.
if there are no equalities of the form
\begin{equation}
\label{e:lc}
v_i=\bigoplus_{j\ne i}\alpha_j v_j.
\end{equation}
The minimal generating sets will be also called \emph{bases}. 

The following crucial result is due to Moller~\cite{Mol88} and
to Wagneur~\cite{Wag88,Wag91}
(it is also contained in~\cite{DevSt} and~\cite{Gau92}).

\begin{proposition}
\label{p:Moller}
If $\{u_1,\ldots,u_s\}$ and $\{v_1,\ldots,v_t\}$ are two bases
of a max-plus space, then $s=t$ and there is a permutation $\sigma$
and a set of nonzero scalars $\{\alpha_1,\ldots,\alpha_s\}$
such that $u_i=\alpha_i v_{\sigma(i)}$ for all
$i=1,\ldots,s$.
\end{proposition}

Prop.~\ref{p:Moller} means that if we have found a finite base for a max-plus
space, then it is in some sense unique: we can only multiply the vectors of the
base by nonzero scalars. For more information on max-plus bases we refer
the reader to~\cite{CunBut04}.

The max-plus matrix algebra is 
formally analogous to the conventional matrix
algebra (minus subtraction and plus idempotency): 
$(A\oplus B)_{ij}=A_{ij}\oplus B_{ij}$ and
$(A\odot B)_{ij}=\bigoplus_k A_{ik}\odot B_{kj}$.

Let us introduce two important characteristics of 
%%%AS
max-plus matrices.
The first characteristic deals with the cyclic permutations.
Let $A$ be an $n\times n$ max-plus matrix. 
%%%AS
Here and below $N$ will stand for the $n$-element set $\{1,\ldots, n\}$.
Denote by $C_n$
the set of all cyclic permutations $\tau$ that act on the subsets
%%%AS
of the set $N$.
%%%AS
For $\tau \in C_n$ denote 
by $K(\tau)$ the subset on which $\tau$ acts and by $\mid K(\tau)\mid$ the
number of elements in this subset. Then
\begin{equation}
\label{e:mcm}
\lambda(A)=\bigoplus_{\tau\in C_n}(\bigodot_{i\in K(\tau)} 
A_{i\tau(i)})^{\frac{1}{\mid K(\tau)\mid}}
\end{equation}
is the \emph{maximal cycle mean} of the matrix $A$
(the notation $\lambda(A)$ for the maximal cycle mean of $A$ will be used
throughout the paper). The summand 
$(\bigodot_{i\in K(\tau)} A_{i\tau(i)})^{\frac{1}{\mid K(\tau)\mid}}$ is called
the \emph{cycle mean} of the cyclic permutation $\tau$. 

The cyclic permutation whose cycle mean
is maximal will be called \emph{critical}.  

The second characteristic deals with the permutations of $N$. 
For any square $n\times n$
max-plus matrix $A$, its \emph{permanent} is defined as
\begin{equation}
\label{e:perm}
\per(A)=\bigoplus_{\sigma\in S_n}\bigodot_{i=1}^n A_{i\sigma(i)}.
\end{equation}
Here $S_n$ is the group of all permutations of $N$.

The summand $\bigodot_{i=1}^n A_{i\sigma(i)}$ is called the \emph{weight of
the permutation} $\sigma$, so the max-plus permanent is the maximal
weight of all permutations of $N$. 
The permanent of $A$ is said to be \emph{strong} (following~\cite{But03}), 
if $A$ has only one \emph{maximal permutation}, i.e., only one permutation with 
maximal weight. 

In the papers~\cite{DSS} and~\cite{DevSt} the max-plus permanent is
called the `tropical determinant'. 
To some extent the max-plus permanent can overtake the role that
the usual determinant plays, 
see~\cite{But03} and the just mentioned papers 
for details. 
We also refer the reader to \cite{BCOQ} and \cite{Gau92},
for a symmetrized version of max-algebra, which admits
subtraction and determinants.

An $n\times n$ max-plus matrix $A$ is called
\emph{invertible}, if there is another $n\times n$ max-plus 
matrix $B$ such that the products $AB$ and $BA$ are both equal to the max-plus identity
matrix $I$. Prop.~\ref{p:Moller} implies that $A$ is
invertible iff there is a permutation $\sigma$ and a set of
nonzero scalars $\alpha_1,\ldots,\alpha_n$ such that
\[
A_{ij}=
\begin{cases}
\alpha_i,  &\text{if $j=\sigma(i)$;}\\
\textbf{0}, &\text{otherwise.}
\end{cases}
\]

So the class of invertible matrices in max-algebra is very small.   
But it makes sense to calculate the series
\begin{equation}
\label{e:closeries}
A^*=I\oplus A\oplus A^2\oplus\ldots,
\end{equation}
where $I$ is the max-plus identity matrix. 
If the sum of this series exists, then it is called the 
\emph{max-plus algebraic closure} of $A$. 
It is an obvious
analogue of $(I-A)^{-1}$. 

%%%AS
Powers  
of max-plus matrices and 
%%%AS
max-plus closures play important role
in optimization on graphs.
Indeed, if we associate an $n$-node graph 
with an $n\times n$ matrix $A$ and let $A_{ij}$ be the weight of the
edge $(i,j)$ of this graph, then the entry $(A^m)_{ij}$ of the matrix
$A^m$ will represent
the maximal weight of paths of length $m$ running from $i$ to $j$. 
The entry $(A^m)_{ii}$ is the maximal weight of all cyclic paths, i.e.
\emph{cycles},
of length $m$ that traverse $i$. In other words, it is the maximal weight
of the cyclic permutations $\tau$ such that $i\in K(\tau)$ and
$\mid K(\tau)\mid=m$.
Analogously, $(A^*)_{ij}$, for $i\ne j$, is the maximal weight of all paths 
of any length running from $i$ to $j$. The path from $i$ to $j$ 
whose weight equals $(A^*)_{ij}$ is called \emph{optimal}. 
%The diagonal entries of $A^*$ 
%do not carry much information.

The following proposition, due to Carr\'{e} \cite{Carre71}, 
solves the problem of existence of the closure.

\begin{proposition}
\label{p:cloexist}
The closure of the matrix $A$ exists if and only if $\lambda(A)\le\textbf 1$.
\end{proposition}

%%%AS
Max-plus  
closures enjoy the property
\begin{equation}
\label{e:multid}
(A^*)^2=A^*.
\end{equation}
Otherwise stated, the inequality $A^*_{ij} A^*_{jk}\le A^*_{ik}$ holds
for all $i$,$j$, and $k$.

As $A^*$ is an equivalent of $(I-A)^{-1}$, 
some algorithms of linear algebra can be adapted
to calculate closures (that is, all optimal paths on a graph) in 
max-algebra, see \cite{Carre71}, \cite{LitMaslova00}, and \cite{Rote85}. 

Another important tool is the max-plus spectral theory. 
For $A$, an $n\times n$ max-plus square matrix, the \emph{max-plus
spectral problem} consists in finding an $n$-element vector $x$, such that
not all of its components are $\textbf 0$, 
and a scalar $\lambda$ such that
\begin{equation}
\label{e:spectral}
Ax=\lambda x.
\end{equation}

The scalar $\lambda$ is a \emph{max-plus eigenvalue}, and the vector $x$ is 
a \emph{max-plus eigenvector}. 
The set of max-plus eigenvectors associated with $\lambda$ is closed
under addition and multiplication by any nonzero scalar. Therefore
it is called the \emph{max-plus eigenspace} (associated with $\lambda$). 
The eigenspace
associated with the maximal eigenvalue of $A$ will be denoted
by $\eig(A)$, and the space generated by the columns of $A$ will
be denoted by $\spann(A)$. 

One of the main results on max-plus spectral theory
is the following.

\begin{proposition}
\label{p:mcm}
The maximal eigenvalue of any max-plus square matrix 
is equal to its maximal cycle mean.  
\end{proposition}

The eigenspace associated with the maximal eigenvalue is easy to describe.

\begin{proposition}
\label{p:eigenspace}
Let $A$ be an $n\times n$ square max-plus matrix such that
$\lambda(A)=\textbf 1$. Then 
\begin{itemize}
\item[1)]
$\eig(A)$ is generated by the columns $A_{\tchk j}^*$
of $A^*$ such that $AA^*_{\tchk j}=A^*_{\tchk j}$;
\item[2)]  $AA^*_{\tchk j}=A^*_{\tchk j}$ iff there is a 
critical cyclic permutation $\tau$
such that $j\in K(\tau)$.
\end{itemize}
\end{proposition}

Columns of~$A^*$ corresponding to vertices of the same cycle are
proportional (in the max-plus sense) to each other.
\begin{proposition}
\label{p:proportia}
Let $A$ be an $n\times n$ max-plus matrix such that 
$\lambda(A)=\textbf 1$ and let $\tau$ be a critical cyclic permutation. 
Then for any
$l\in N$ and $i\in K(\tau)$ we have
\begin{equation}
\begin{array}{l}
A_{il}^* = A_{i\tau(i)} A_{\tau(i)l}^*,\\
A_{li}^* = A_{l\tau^{-1}(i)}^* A_{\tau^{-1}(i)i}.
\end{array}
\end{equation}
\end{proposition}

If the graph associated with $A$ is strongly connected (i.e. for
any $i$ and $j$ there is a path with nonzero weight running from $i$ to $j$),
then $A$ is said to be \emph{irreducible}. In this case the maximal
cycle mean of $A$ is known to be its only eigenvalue. 
In the reducible case 
this eigenvalue need not be unique.  
If another eigenvalue exists, then it is the (only) eigenvalue
of some maximal irreducible submatrix of $A$. 
However, not all
maximal irreducible submatrices of $A$ yield 
%%%AS 
eigenvalues for $A$.

For more details on max-plus spectral theory, as well as for
the proofs of the above propositions, we refer the reader to
\cite{BCOQ}, \cite{CG79}, \cite{Gau92}, and~\cite{Zimm81}.

\section{Definite eigenspaces and definite closures}
\label{s:defclos}
We begin this section with the following important proposition.
The proof makes use of the uniqueness of the base (Prop.~\ref{p:Moller}).
\begin{proposition}
\label{p:*=*}
Let $A$ and $B$ be square max-plus matrices such that $\lambda(A)\le\textbf 1$,
and $\lambda(B)\le\textbf 1$. Then $\spann(A^*)=\spann(B^*)$ if and only if
$A^*=B^*$.
\end{proposition}

\textbf{Proof.}
First let us prove that whenever
\begin{equation}
\label{e:lcomb}
A^*_{\tchk i}=\bigoplus_{j\ne i} \alpha_j A^*_{\tchk j},
\end{equation}
the column~$A^*_{\tchk i}$ is proportional to~$A^*_{\tchk j}$ for
some~$j\neq i$.

%%%AS
Suppose (\ref{e:lcomb}) holds.  Then  
there is an~$l$ such 
that $A_{ii}^*=\textbf{1}=\alpha_l A^*_{il}$ and that
$A^*_{li}\ge\alpha_l$. Combining this we obtain that 
$A_{il}^* A_{li}^*\ge\textbf 1$, hence $A_{il}^* A_{li}^*=\textbf 1$
(otherwise there is a cycle whose weight exceeds $\textbf 1$). This
means that there is a critical cycle with weight $\textbf 1$
traversing $i$ and $l$, and due to Prop.~\ref{p:proportia} 
$A^*_{\tchk i}=\alpha_l A^*_{\tchk l}$. We conclude that no
column of $A^*$ can be expressed as linear combination~(\ref{e:lcomb})
without being proportional to some of the columns involved in this combination.

Further let $\{\,A^*_{\tchk r_1},\ldots , A^*_{\tchk r_k}\,\}$ be the
base of $\spann(A^*)=\spann(B^*)$. If we use
%%%AS 
columns of $B^*$
to form the base, then, due
to Prop.~\ref{p:Moller}, it must be of the form
$\{\,B^*_{\tchk s_1},\ldots , B^*_{\tchk s_k}\,\}$, so that
$B^*_{\tchk s_i}=\alpha_i A^*_{\tchk r_{\sigma(i)}}$ for $i=1,\ldots ,k$
and some nonzero $\alpha_i$.
All remaining columns of $A^*$ (of $B^*$) are   
proportional to 
%%%AS 
base columns of $A^*$ (of $B^*$).
So every column of $A^*$ is proportional to some column of $B^*$,
and vice versa. This implies that the rows $A^*_{i\tchk}$ and
$A^*_{j\tchk}$ are proportional iff the rows $B^*_{i\tchk}$
and $B^*_{j\tchk}$ are proportional. 

Let the columns $A^*_{\tchk i}$ and $A^*_{\tchk j}$ be proportional.
Then $A_{ij}^* A_{ji}^*=\textbf 1$, hence there is a critical cycle
containing $i$ and $j$. Due to Prop.~\ref{p:proportia} 
the rows $A^*_{i\tchk}$ and $A^*_{j\tchk}$
are proportional, so are the rows $B^*_{i\tchk}$ and $B^*_{j\tchk}$,
and, again due to Prop.~\ref{p:proportia}, so are the columns $B^*_{\tchk i}$
and $B^*_{\tchk j}$. We conclude that the columns $A^*_{\tchk i}$
and $A^*_{\tchk j}$ are proportional iff so are the columns $B^*_{\tchk i}$
and $B^*_{\tchk j}$.

Now it is clear that $\{\,A^*_{\tchk r_1},\ldots,A^*_{\tchk r_k}\,\}$
is the base of $\spann(A^*)$ iff  
$\{\,B^*_{\tchk r_1},\ldots,B^*_{\tchk r_k}\,\}$
is also the base, so that $B^*_{\tchk r_i}=\alpha_i 
A^*_{\tchk r_{\sigma(i)}}$. We can assume w.l.o.g.\ that $r_i=i$.
 Consider the decomposition of $\sigma$
into  
cyclic permutations and let $\tau$ be one of them. Then
$B^*_{\tau(i)i}=\alpha_i$ and $A^*_{i\tau(i)}=\alpha_i^{-1}$ for
all $i\in K(\tau)$. If $\bigodot_{i\in K(\tau)}\alpha_i>\textbf 1$, then
$B^*$ has a cycle with weight greater than $\textbf 1$, and if 
$\bigodot_{i\in K(\tau)}\alpha_i<\textbf 1$, then so does $A^*$. The only remaining
possibility $\bigodot_{i\in K(\tau)} \alpha_i =\textbf 1$ implies
that all columns of $A^*$ and $B^*$ with indices belonging to $K(\tau)$
are proportional. Then $K(\tau)$ must be a singleton for any $\tau$,
otherwise the minimality of the bases is violated. This implies that
$\sigma$ is the identity permutation, and $A^*_{\tchk i}=B^*_{\tchk i}$
for any $i=1,\ldots,k$. 
Taking into account that there is a one-to-one correspondence between
the sets of proportional columns of $A^*$ and $B^*$, and that all
columns of $A^*$ and $B^*$ are proportional to some base columns, we
conclude that $A^*=B^*$.\hfill$\square$  

Now consider 
%%%AS 
matrices with maximal cycle mean equal to $\textbf 1$
and with all diagonal entries equal to $\textbf 1$. 
Following~\cite{But03}, we call such matrices \emph{definite}.
The following
proposition contains some simple facts on 
%%%AS 
eigenspaces of such
matrices. The third statement is an easy consequence of~\cite{Gau92}, Ch.~IV,
Th.~2.2.4, and its proof is 
%%%AS
recalled here  
for convenience of the reader.
We consider the general reducible case, 
%%%AS
in which 
the eigenvectors
may have zero entries. For any $y\in \rn_{\max}^n$, the index set $K$
such that $y_i\ne\textbf 0$ 
%%%AS
iff  
$i\in K$ is called the \emph{support}
of $y$ and 
%%%AS
is
denoted by $\supp(y)$.

\begin{proposition}
\label{p:defmat}
If $A$ is a definite matrix, then
\begin{itemize}
\item[1)] it has 
%%%AS
a
unique eigenvalue equal to $\textbf 1$;
\item[2)] $\eig(A)=\spann(A^*)$;
\item[3)] an eigenvector with the support $K\subset N$ exists iff
$A_{ij}=\textbf 0$ for all $i\in N\bez K$ and $j\in K$.
\end{itemize}
\end{proposition}

\textbf{Proof.} 1) Any eigenvalue of $A$ is the maximal cycle mean
of some of its submatrices, and the maximal cycle mean of any
submatrix of $A$ is equal to $\textbf 1$.

2) 
%%%AS 
follows from Prop.~\ref{p:eigenspace}.

3) If an eigenvector $x$ such that $\supp(x)=K$ exists, then $A_{ij} x_j=\textbf 0$
for all $i\in N\bez K$ and all $j\in K$. Hence $A_{ij}=\textbf 0$ for
all such $i$ and $j$.

Conversely, if the set $K$ satisfies the condition, we look for
eigenvectors $x$ such that 
$\supp(x)=K$. We may reduce the system $Ax=x$ to the
system $A_{KK} x_K=x_K$, where $A_{KK}$ is a submatrix of $A$ standing
on the rows and columns with indices belonging to $K$, and $x_K$ is a vector
with $\mid K\mid$ nonzero components. The space $\eig(A_{KK})$ is generated
by the columns $(A_{KK})^*_{\tchk j}$, where $j\in K$. Taking any
combination of all these generators with all coefficients not equal to
$\textbf 0$ we obtain an eigenvector with the support $K$.\hfill$\square$

In the case of definite matrices we have one more implication
of the uniqueness of the base. The proof is similar to that of
Prop.~\ref{p:*=*}.

\begin{proposition}
\label{p:A=*}
If $A$ is definite and $\spann(A)=\eig(A)$, then $A=A^*$.
\end{proposition}

\textbf{Proof.} Let $\{\,A_{\tchk s_1},\ldots A_{\tchk s_k}\,\}$ be the
base of $\spann(A)=\eig(A)=\spann(A^*)$. Due to Prop.~\ref{p:Moller}, 
if we use 
%%%AS 
columns of
$A^*$ to form the base, then it
must be of the form  $\{\,A_{\tchk t_1}^*,\ldots A_{\tchk t_k}^*\,\}$ so
that $A_{\tchk s_i}=\alpha_i A^*_{\tchk t_i}$ for $i=1,\ldots,k$ and some
nonzero $\alpha_i$. More precisely, $\alpha_i=A_{t_is_i}$, and this
implies $A_{t_is_i}A^*_{s_it_i}=\textbf 1$. Then there is a critical cycle
traversing $s_i$ and $t_i$, so $A^*_{\tchk s_i}=A_{t_is_i} A_{\tchk t_i}^*$
according to Prop.~\ref{p:proportia}. So $A_{\tchk s_i}=A_{\tchk s_i}^*$
for all $i=1,\ldots,k$.

Now assume that there are columns $A_{\tchk j}$ and scalars $\alpha_i$ 
such that
\begin{equation}
\label{e:alcomb}
A_{\tchk j}=\bigoplus_{i=1}^k \alpha_i A_{\tchk s_i},
\end{equation}
and that no column of the base is proportional to $A_{\tchk j}$. It
follows from~(\ref{e:alcomb}) that there is an index $m$ such that
$A_{jj}=\textbf{1}=\alpha_m A_{js_m}$. The columns $A_{\tchk j}$ and
$A_{\tchk s_m}$ are not proportional, hence there is an $l$ such that
$A_{lj}>\alpha_m A_{ls_m}$. This implies $A_{lj}A_{js_m}>A_{ls_m}$
and $A^*_{ls_m}>A_{ls_m}$. This is a contradiction, since we have
proved that $A^*_{\tchk s_i}=A_{\tchk s_i}$ for any $i=1,\ldots,k$. 
So any column of $A$ is proportional to some column of the base.
But if $A_{\tchk i}$ and $A_{\tchk j}$ are proportional, then 
$A_{ij}A_{ji}=\textbf 1$, hence $A^*_{\tchk i}$ and $A_{\tchk j}^*$
are also proportional with the same coefficient. This implies $A=A^*$.
\hfill$\square$

Let us now introduce a \emph{definite form} of a matrix.
Consider an $n\times n$ max-plus matrix $A$ that has nonzero permanent,
i.e., at least one permutation whose weight is not equal to $\textbf 0$.
Let $\sigma$ be one of the maximal permutations of $A$.
If this permutation is not the identity permutation, then we turn this
permutation into the identity permutation by rearranging the columns of $A$.
Then we divide (in the max-plus sense) 
all columns by the corresponding diagonal entries, thus
obtaining 
%%%AS
a  
matrix $A'$ with entries
\begin{equation}
\label{e:deform}
A'_{ij}=A_{i\sigma(j)} A_{j\sigma(j)}^{-1}.
\end{equation}
Obviously, passing to a definite form of a matrix does not alter
its span: $\spann(A') = \spann(A)$.

It is clear that $A'$ is definite, hence we call it the \emph{definite form} 
of $A$ corresponding to the permutation $\sigma$. For example, if
\[
A=
\begin{pmatrix}
1&5&-\infty\\
2&1&7\\
6&-\infty&2
\end{pmatrix},
\]
then
\[
A'=
\begin{pmatrix}
0&-\infty&-5\\
-4&0&-4\\
-\infty& -5 &0
\end{pmatrix}
\]
is the definite form of $A$ corresponding to $\sigma=(2 3 1)$.
(Here and in the sequel we prefer not to use the symbols $\textbf 0$ and
$\textbf 1$ in numerical examples.) 

However, in general there are many maximal permutations and many
definite forms corresponding to them. The fact that we want to prove
is that closures of all definite forms coincide. It is convenient to
pose the problem as follows. Let $A$ be definite, 
let it have maximal permutations different
from the identity permutation, and let $\sigma$ be one of them. Let
$A'$ be the definite form of $A$ corresponding to $\sigma$, then its
entries are defined according to~(\ref{e:deform}). $A'$ also has maximal
permutations different from the identity permutation, and $\sigma^{-1}$
is one of them, since
\begin{equation}
\label{e:sigma-1}
A'_{i\sigma^{-1}(i)}=A_{\sigma^{-1}(i)i}^{-1}.
\end{equation}

First we prove the following proposition.

\begin{proposition}
\label{p:eig=eig}
Let $A$ be a definite matrix and let $\sigma$ be one of its
maximal permutations. For $A'$ the definite form of $A$ corresponding to
$\sigma$, $\eig(A')=\eig(A)$.
\end{proposition}

\textbf{Proof.} 
Consider the decomposition of $\sigma$ into cyclic permutations.
Let $r$ be the number of these permutations and 
denote these permutations by $\tau_l$ (where $l=1,\ldots,r$).  Denote by 
$K(\tau_l)$ the index set on which $\tau_l$ acts. 
The sets $K(\tau_l)$ are pairwise disjoint and
$\bigjoin_{l=1}^r K(\tau_l) = N$. 
%The inverse permutation 
%$\sigma^{-1}$ is decomposed into the
%inverse cyclic permutations $\tau^{-1}_l$, and $K(\tau^{-1}_l)=K(\tau_l)$.

We prove the inclusion $\eig(A)\subset\eig(A')$.
Let $y$ be an eigenvector of $A$ with the support $M\subset N$. Due
to the third statement of Prop.~\ref{p:defmat} the set $M$ must be of the
form $\bigjoin_{l\in L} K(\tau_l)$ for some 
$L\subset\{1,\ldots,r\}$, and we have
$A_{ij}=\textbf 0$ for all $i\in N\bez M$ and $j\in M$.

Now note that every cyclic permutation $\tau_l$ 
is critical, hence for any $\tau_l$ such that $K_l\in M$ we have
\begin{equation}
\label{e:big=1}
\bigodot_{i\in K(\tau_l)} A_{i\tau_l(i)} y_{\tau_l(i)} y_i^{-1}=\textbf{1}.
\end{equation}

On the other hand we have $A_{i\tau_l(i)} y_{\tau_l(i)}\le y_i$.
This inequality is an equality, otherwise the violation of~(\ref{e:big=1})
occurs. Thus we obtain
\begin{equation}
\label{e:cycle=}
A_{i\sigma(i)} y_{\sigma(i)}=y_i
\end{equation}
for all $i\in M$.

Now note that, since $\sigma(j)\in M$ for all $j\in M$, we have that
$A'_{ij}=\textbf 0$ for all $i\in N\bez M$ and $j\in M$, 
the same as for $A_{ij}$. 
Therefore it suffices to prove that if $A_{MM} y_M=y_M$, then $A'_{MM} y_M=y_M$.
In other words, we want to show that $A'_{ij} y_j\le y_i$ for all $i,j\in M$.
But due to~(\ref{e:deform})
and~(\ref{e:cycle=})
\[
\begin{array}{l@{{}={}}l}
A'_{ij}y_j &A_{i\sigma(j)} A_{j\sigma(j)}^{-1} y_j=\\[1.5ex]
& A_{i\sigma(j)} y_{\sigma(j)} y_{\sigma(j)}^{-1} A_{j\sigma(j)}^{-1} y_j=\\[1.5ex]
& A_{i\sigma(j)} y_{\sigma(j)}\le y_i,
\end{array}
\]
and the proof of the inclusion $\eig(A)\subset\eig(A')$ is complete. The 
opposite inclusion is proved analogously, after passing to the definite
form associated with $\sigma$ (the matrices $A$ and $A'$ will interchange).  
\hfill$\square$

The fact that we want to prove is now an immediate consequence of
Propositions~\ref{p:*=*} and~\ref{p:eig=eig}.

\begin{proposition}
\label{p:defclos}
Closures of all definite forms of any matrix with nonzero permanent
coincide.
\end{proposition}

We have just proved that closures of all definite forms coincide.
%%%AS 
Now we can define the \emph{definite closure} of any $n\times n$
max-plus matrix with nonzero permanent to be the closure of any
of its definite forms.
%%%AS
Due to the second
statement of Prop.~\ref{p:defmat}, for any
definite form $A'$ of $A$ we have $\eig(A')=\eig((A')^*)=\spann((A')^*)$.
So the eigenspace of the definite closure of $A$ coincides with the eigenspace
of any of its definite forms 
%%%AS
(they are all the same),
%%%AS
and is generated
by  
columns of the definite closure. The eigenspace of definite closure
will be called the \emph{definite eigenspace}.
 
The third statement  
of Prop.~\ref{p:defmat} suggests
that, if we want to work with the space of 
eigenvectors of definite closure 
%%%AS 
of non-full 
%%%AS
support, then we must confine ourselves to the corresponding
submatrix. 

Further we always assume that the eigenvectors considered
have full support, i.e., that we study the eigenvectors with certain
support and have passed to the corresponding submatrix.

Let us show that the definite eigenspace can be described
by some system of inequalities.

\begin{proposition}
\label{p:syseq}
Let $A$ be a definite matrix. Then its eigenspace (and the eigenspace
of its closure) is the set 
$X=\{\, x\mid A_{ij}\le x_ix_j^{-1},\, i\ne j,\, A_{ij}=A^*_{ij}\,\}$ 
\end{proposition}  
\textbf{Proof.} Let $x$ be an eigenvector of $A$ corresponding to
its maximal eigenvalue $\textbf 1$. It satisfies the system $Ax=x$,
therefore it satisfies all inequalities of the form $A_{ij}\le x_ix_j^{-1}$.
Hence $x\in X$.

Conversely, if all inequalities $A_{ij}\le x_ix_j^{-1}$, for $i\ne j$
and $A_{ij}=A_{ij}^*$, are satisfied, then absolutely all inequalities
$A_{ij}\le x_ix_j^{-1}$ are satisfied. Indeed, if $A_{ij}<A_{ij}^*$, then the
optimal path from $i$ to $j$ is not the edge $(i,j)$, but it traverses
other nodes, say, $i_1,\ldots,i_k$. Then
$A_{ij}^*=A_{ii_1}\ldots A_{i_kj}$ where $A_{ii_1}=A_{ii_1}^*,\ldots,
A_{i_kj}=A_{i_kj}^*$ (all edges $(i,i_1),\ldots,(i_k,j)$ are optimal paths).
The inequality $A_{ij}\le x_ix_j^{-1}$ is now an easy consequence
of the inequalities $A_{ii_1}\le x_ix_{i_1}^{-1},\ldots,\
A_{i_kj}\le x_{i_k}x_j^{-1}$ that are satisfied. So all inequalities
$A_{ij}\le x_ix_j^{-1}$ are satisfied and this implies $Ax=x$.
\hfill$\square$

An inverse problem can also be posed. Suppose that there is a system
of inequalities $\{a_{ij}\le x_ix_j^{-1}\}$, with at most one inequality
per each pair $(i,j)$, and the set of vectors with
full support defined by this system is not empty.
Is it a full-support subspace of an eigenspace of a definite matrix?

To answer this question, consider the matrix $A$ whose entries $A_{ij}$ are
equal either to $\textbf 1$, if $i=j$, or to $a_{ij}$, 
if there is an inequality of
the form $a_{ij}\le x_ix_j^{-1}$, or to $\textbf 0$, if there is no such
inequality. Then we have the following proposition.

\begin{proposition}
\label{p:invprob}
The set $X=\{\, x\mid A_{ij}\le x_ix_j^{-1}\,\}$ is nonempty if and only if
$A$ is definite.
\end{proposition} 
\textbf{Proof.} If $A$ is definite then the set of its eigenvectors
associated with the eigenvalue $\textbf 1$ is nonempty and, due
to Prop.~\ref{p:syseq}, it is precisely the set $X$ (some of inequalities
being redundant).

Conversely, let $X$ be nonempty, then there exists $x\in X$. Take an
arbitrary cyclic permutation $\tau$. As a consequence of all inequalities 
$A_{i\tau(i)}\le x_ix_{\tau(i)}^{-1}$, where $i\in K(\tau)$, we obtain that
$\bigodot_{i\in K(\tau)} A_{i\tau(i)}\le\textbf 1$. Hence all cycle means
of $A$ are not greater than $\textbf 1$, and $A$ is definite.\hfill$\square$
 
Now consider the following application.

Let $V$ be an $m\times n$ max-plus matrix with at least one
nonzero entry in each column,
and let $y\in\rn_{\max}^m$ have full support. Then, following~\cite{DevSt},
we can define the \emph{combinatorial type} of $y$ with respect to $V$.
It is an $m$-tuple $S$ of subsets $S_1,\ldots,S_m$ of $\{1,\ldots,n\}$,
such that $i\in S_j$ whenever 
$\bigoplus_{k=1}^m V_{ki}y_k^{-1}=V_{ji}y_j^{-1}$. It is
proved in~\cite{DevSt} that the collection of the sets
\begin{equation}
\label{e:cells}
X_S=\{\,y\mid V_{ki}V_{ji}^{-1}\le y_ky_j^{-1}\,\text{for $j,k=1,\ldots,m$ and 
$i\in S_j$}\,\}
\end{equation}
defines a cellular decomposition of the full-support subspace of 
$\rn_{\max}^m$. 

Consider an $m\times m$ matrix $V^S$ whose columns are defined by
\begin{equation}
\label{e:vnonempty}
V^S_{\tchk j}=\bigoplus_{i\in S_j} V_{ji}^{-1}V_{\tchk i},
\end{equation}
if $S_j$ is not empty, and by
\begin{equation}
\label{e:vempty}
V^S_{ij}=
\begin{cases}
\textbf{1}, &\text{if $i=j$;}\\
\textbf{0}, &\text{if $i\ne j$.}
\end{cases}
\end{equation}
if $S_j$ is empty.
Then it is clear that
\begin{equation}
\label{e:xsvs}
X_S=\{\,y\mid V_{ij}^S\le y_iy_j^{-1}\,\text{for $i,j=1,\ldots,m$}\,\},
\end{equation}
and we immediately have the following proposition.

\begin{proposition}
\label{p:Sturmfels}
\begin{itemize}
\item[1)] The cell $X_S$ 
%%%AS 
is nonempty  
%%%AS
if and only if $\lambda(V^S)\le\textbf 1$;
\item[2)] If $\lambda(V^S)\le\textbf 1$, then 
$X_S=\{\,y\in\spann((V^S)^*)\mid \supp(y)=N\,\}$;
\item[3)] If $\lambda(V^S)\le\textbf 1$, then 
$X_S=\{\, x\mid V^S_{ij}\le x_ix_j^{-1},\, i\ne j,\, V^S_{ij}=(V^S)^*_{ij}\,\}$.
\end{itemize}
\end{proposition}

We close this section with three examples.
\vskip1cm
\emph{Example 1} Consider the matrix
\[
A=
\begin{pmatrix}
1&3&0\\
2&0&0\\
0&-1&-5
\end{pmatrix}.
\]

The only maximal permutation of $A$ is $(231)$, and the only
definite form is
\[
A'=
\begin{pmatrix}
0&0&1\\
-3&0&2\\
-4&-5&0
\end{pmatrix}.
\]

The definite closure of $A$ is
\[
(A')^*=
\begin{pmatrix}
0&0&2\\
-2&0&2\\
-4&-4&0
\end{pmatrix}.
\]

Fig.~1 displays the cross section by $z=0$ of
$\spann(A)$ (left) and $\spann((A')^*)=\eig(A')$ (right).

%\begin{figure}[h]  
%\centering
%\vskip-5cm
%\includegraphics[width=10cm]{dcfora}  
%\vskip-5cm
%\caption{The definite closure operation for $A$}
%\end{figure}

\begin{figure}[htbp]
\label{f:stper}
%\vskip0.5cm
\hskip1cm
\epsfysize=4cm
\epsffile{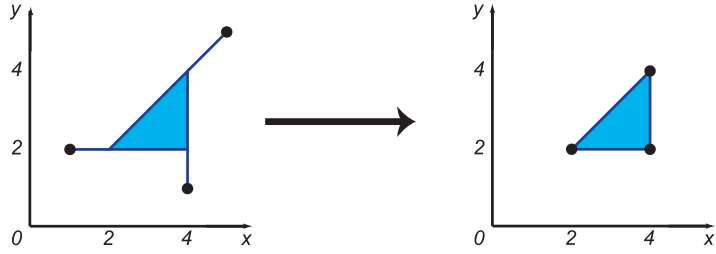}
%\vskip0.5cm
\caption{The definite closure operation for $A$}
\end{figure}

Note that $\eig(A')$ is the set 
\[
\{(x,y,z)\mid\;xy^{-1}\ge 0,\; yz^{-1}\ge 2,\; zx^{-1}\ge -4\;\},
\]
in accordance
with Prop.~\ref{p:syseq}.
\vskip1cm
\emph{Example 2} Consider the matrix
\[
B=
\begin{pmatrix}
2&0&2\\
1&1&3\\
0&-3&-2
\end{pmatrix}
\]

It has two maximal permutations: $(13)(2)$ and $(231)$. Therefore
it has two definite forms, namely
\[
B'=
\begin{pmatrix}
0&-1&2\\
1&0&1\\
-4&-4&0
\end{pmatrix}
\]

and
\[
B''=
\begin{pmatrix}
0&-1&2\\
1&0&1\\
-3&-5&0
\end{pmatrix}.
\]

But, in accordance with Prop.~\ref{p:defclos}, the definite closure
\[
(B')^*=(B'')^*=
\begin{pmatrix}
0&-1&2\\
1&0&3\\
-3&-4&0
\end{pmatrix}
\]
is unique.

Fig.~2 displays the cross section by $z=0$ of
$\spann(B)$ (left) and $\spann((B')^*)=\eig(B')$ (right).

%\begin{figure}[h]  
%\centering
%\vskip-5cm
%\includegraphics[width=10cm]{dcforb}  
%\vskip-5cm
%\caption{The definite closure operation for $B$}
%\end{figure}

\begin{figure}[htbp]
\label{f:wper}
%\vskip0.5cm
\hskip1cm
\epsfysize=4cm
\epsffile{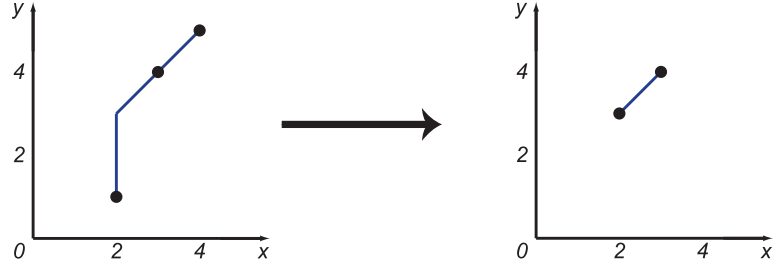}
%\vskip0.5cm
\caption{The definite closure operation for $B$}
\end{figure}

In accordance with Prop.~\ref{p:syseq}, the space $\eig(B')$ is the set 
\[
\{\,(x,y,z)\mid yx^{-1}\ge 1,\, xy^{-1}\ge -1,\, xz^{-1}\ge 2,\,
zy^{-1}\ge -4\,\},
\] 
or, equivalently, the set
\[
\{\,(x,y,z)\mid yx^{-1}\ge 1,\, xy^{-1}\ge -1,\, xz^{-1}\ge 2,\,
zx^{-1}\ge -3\,\}.
\]
The first system of inequalities corresponds to the definite form $B'$
and the second one corresponds to $B''$.

Comparing Fig.~1 with Fig.~2 we see that $\eig(A')$
has `interior' whereas $\eig(B')$ does not have `interior' (for the
exact meaning of the term `interior' see Sect.~\ref{s:hdist} below). 
As a consequence of~\cite{But00}, Th.~4.2, or
~\cite{DSS}, Th.~4.2, one can obtain that 
$\eig(A')$, for $A'$ a definite form of $A$,  
has `interior' if and only if $A$ (or equivalently $A'$) has
strong permanent. 
This fact will be revisited in Prop.~\ref{p:but} of this paper. 

\vskip1cm

\emph{Example 3} Consider the matrix
\[
V=
\begin{pmatrix}
1&4&6&7\\
4&1&5&8\\
0&0&0&0
\end{pmatrix}.
\]
Let $S=(\,\{2,3\},\{4\},\{1\}\,)$, $P=(\,\emptyset,\{4\},\{1,2,3\}\,)$,
$U=(\,\{3\},\{1,3,4\},\{2\}\,)$,
and $W=(\,\{1,4\},\{2\},\{3\}\,)$ be four combinatorial types.
Do they exist in the cellular decomposition? If they do, what
vectors generate the respective cells $X_S$, $X_P$, $X_U$, and $X_W$?

For $S$:
\[
V^S=
\begin{pmatrix}
0&-1&1\\
-1&0&4\\
-4&-8&0
\end{pmatrix},\quad
(V^S)^*=
\begin{pmatrix}
0&-1&3\\
0&0&4\\
-4&-5&0
\end{pmatrix}.
\]

Hence $X_S$ exists and is generated by $[4\;4\;0]^T$,
$[4\; 5\; 0]^T$, and $[3\;4\;0]^T$.

For $P$:
\[
V^P=
\begin{pmatrix}
0&-1&6\\
-\infty&0&5\\
-\infty&-8&0
\end{pmatrix}
=(V^P)^*.
\] 

Hence $X_P$ is generated by $[0 -\infty -\infty]^T$,
$[7\; 8\; 0]^T$ and $[6\; 5\; 0]^T$. However,
$[0 -\infty -\infty]^T$ does not have full support
and does not belong to $X_P$.

For $U$:

\[
V^U=
\begin{pmatrix}
0&1&4\\
-1&0&1\\
-6&-4&0
\end{pmatrix},\quad
(V^U)^*=
\begin{pmatrix}
0&1&4\\
-1&0&3\\
-5&-4&0
\end{pmatrix}.
\]

Hence $X^U$ exists and is generated by $[5\;4\;0]^T$ and $[4\;3\;0]^T$.

For $W$:
\[
V^W=
\begin{pmatrix}
0&3&6\\
3&0&5\\
-1&-1&0
\end{pmatrix}.
\]

We have $\lambda(V^W)=3>\textbf 1$, hence $W$ does not exist in the
cellular decomposition.

Fig.~3 displays $\spann(V)$ (blue), $X_S$ (red), $X_P$ (light grey),
and $X_U$(dark green),
projected onto $z=0$.
The generators of $\spann(V)$ (the 
%%%AS 
larger 
%%%AS 
circles) and the generators
of $X_S$ and $X_U$(the smaller squares) are also shown.

%\begin{figure}[h]  
%\centering
%\vskip-2cm
%\includegraphics[width=8cm]{cells}  
%\vskip-3cm
%\caption{Three cells of a cellular decomposition}
%\end{figure}

\begin{figure}[htbp]
\label{f:cells}
%\vskip0.5cm
\hskip1cm
\epsfysize=5cm
\epsffile{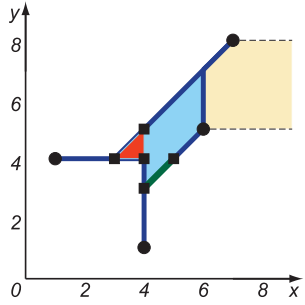}
%\vskip0.5cm
\caption{Three cells of a cellular decomposition}
\end{figure}

\section{Inner structures of definite eigenspaces}
\label{s:hdist}
We have introduced the definite closure
operation and have given an external description of a definite eigenspace
in terms of a system of inequalities. In this
section we give an internal description of definite eigenspace and
measure Hilbert distances between the structures involved in this
description
%%%AS
and the boundary.

For further considerations we need the following notions and
notation. Let $A$ be a definite matrix. 
The sets $X_{ij}=\{\,x\mid A_{ij}=x_ix_j^{-1}\,\}$ will
be called the \emph{supporting planes} of $\eig(A)$, and the
sets $\Gamma_{ij}=X_{ij}\meet \eig(A)$ will be called
the \emph{faces} of $\eig(A)$. The \emph{boundary}, i.e., the union of all faces
of $\eig(A)$ will be denoted by $\Gamma(\eig(A))$. The set of eigenvectors
not belonging to the boundary will be called the \emph{interior} of
$\eig(A)$ and denoted by $\inter(\eig(A))$. 

Also, denote by $\A$ the matrix
obtained from $A$ by replacing the diagonal $\textbf 1$'s by
$\textbf 0$'s, and denote by $A_{\mu}$ the matrix obtained
from $\mu^{-1}{\A}$ by replacing the diagonal $\textbf 0$'s by
$\textbf 1$'s. For example, if
\[
A=
\begin{pmatrix}
0&-4&1\\
1&0&1\\
-5&-7&0
\end{pmatrix}
\]
then
\[
\A=
\begin{pmatrix}
-\infty&-4&1\\
1&-\infty&1\\
-5&-7&-\infty
\end{pmatrix}
\]
and, e.g.,
\[
A_{-1}=
\begin{pmatrix}
0&-3&2\\
2&0&2\\
-4&-6&0
\end{pmatrix}
\]
Also, the maximal cycle mean of $\A$ is $\lambda(\A)=-1.5$, and
\[
A_{\lambda}=
\begin{pmatrix}
0&-2.5&2.5\\
2.5&0&2.5\\
-3.5&-5.5&0
\end{pmatrix}
\]
(we write $A_{\lambda}$ instead of $A_{\lambda(\A)}$ for the sake of
simplicity)

The following crucial proposition can be derived from~\cite{But00}. 
%We prove it
%again for convenience of the reader.

\begin{proposition}
\label{p:but}
Let $A$ be a definite matrix such that 
$\lambda(\A)\ne\textbf 0$.
\begin{itemize}
\item[1)] If $A$ does not have strong permanent then $\eig(A)$ does
not have interior;
\item[2)] If $A$ has strong permanent then $\eig(A)$ has interior, and
\begin{equation}
\label{e:repbut}
\inter(\eig(A))=\bigjoin_{\lambda(\A)\le\mu<\textbf 1} \eig(A_{\mu}).
\end{equation}
\end{itemize}
\end{proposition}
\textbf{Proof.} 1)
If $A$ does not have strong permanent, then it has
a maximal permutation that differs from the identity permutation.
Let $\sigma$ be such permutation, and assume that there is a $y$
belonging to $\inter(\eig(A))$. Then $A_{i\sigma(i)}y_{\sigma(i)}<y_i$
for all $i$. After multiplying all these inequalities and then cancelling
the product $\bigodot_i y_i$
one obtains $\bigodot_i A_{i\sigma(i)} < \textbf 1$. This implies
that $\sigma$ is not maximal, a contradiction.

2)
$y$ belongs to $\inter(\eig(A))$ if and only if 
$\bigoplus_{j\ne i} A_{ij}y_j<y_i$ for all $i$. This takes place
if and only if there is a $\mu<\textbf 1$ such that
$\bigoplus_{j\ne i} A_{ij}y_j\le\mu y_i$, or, equivalently,
$A_{\mu}y=y$. If we take $\mu\ge\lambda(\A)$, then $\textbf 1$ is
the maximal cycle mean of $A_{\mu}$, hence
$\eig(A_{\mu})$ exists. The proof is complete.

\hfill$\square$

If $\lambda(\A)=\textbf 0$, i.e. if the graph associated with $\A$
is acyclic, then representation~(\ref{e:repbut}) is
replaced by the representation
\begin{equation}
\label{e:repbut1}
\inter(\eig(A))=\bigjoin_{\textbf{0}<\mu<\textbf 1} \eig(A_{\mu}).
\end{equation}

In the sequel, we always assume that $A$ is definite and has at
least one off-diagonal entry not equal to $\textbf 0$. 

According to Prop.~\ref{p:syseq}, the eigenspace $\eig(A)$ is the set
$$X=\{\,x\mid A_{ij}\le x_ix_j^{-1},\, i\ne j,\, A_{ij}=A_{ij}^*\,\}.$$

Analogously, the eigenspace $\eig(A_{\mu})$, for any $\mu$ involved
in~(\ref{e:repbut}) or~(\ref{e:repbut1}), is the set

$$X_{\mu}=\{\,x\mid \mu^{-1}A_{ij}\le 
x_ix_j^{-1},\, i\ne j,\, \mu^{-1}A_{ij}=(A_{\mu}^*)_{ij}\,\}.$$

We need the following proposition
mainly for the proof of Prop.~\ref{p:=}.

\begin{proposition}
\label{p:injection}
Let $\mu$ be a scalar such that
$\lambda(\A)\le\mu<\textbf 1$, if $\lambda(\A)>\textbf 0$,
or such that $\textbf{0}<\mu<\textbf 1$, if $\lambda(\A)=\textbf 0$. 
If $(A_{\mu}^*)_{ij}=(A_{\mu})_{ij}$, then $A^*_{ij}=A_{ij}$.
\end{proposition}
\textbf{Proof.} In both cases considered the maximal cycle mean
of $A_{\mu}$ is equal to $\textbf 1$, hence $A_{\mu}^*$ exists.
Let $A^*_{ij}>A_{ij}$, then 
$A^*_{ij}=A_{ii_1}\ldots A_{i_kj}$ for some $i_1,\ldots,i_k$
not equal to $i$ or $j$. Since $\mu<\textbf 1$, we have 
$\mu^{-1}A_{ii_1}\ldots \mu^{-1}A_{i_kj}>\mu^{-1}A_{ij}$,
hence $(A_{\mu}^*)_{ij}>(A_{\mu})_{ij}$.\hfill$\square$

The eigenspaces $\eig(A_{\mu})$ are the inner structures mentioned
above. Now we are going to measure the Hilbert distances between these
inner structures and the boundary $\Gamma(\eig(A))$.

The \emph{Hilbert 
distance} between the two vectors 
$x$ and $y$ both having support $K$ is defined to be
\begin{equation}
\label{e:hdist}
d_H(x,y)=\bigoplus_{i,j\in K} x_ix_j^{-1}y_i^{-1}y_j.
\end{equation}

Note that in~\cite{CGQ04} the Hilbert distance is defined as an inverse
of the quantity $d_H(x,y)$. 
If the supports of $x$ and $y$ differ, then we assume 
the Hilbert distance between
$x$ and $y$ to be infinite.

It can be easily verified (see also~\cite{CGQ04}, Th.~17) 
that the following properties
hold:
\begin{itemize}
\item[1)] $d_H(x,y)\ge\textbf 1$, and $d_H(x,y)=\textbf 1$ iff $x=\lambda y$,
where $\lambda$ is a nonzero scalar;
\item[2)] $d_H(x,y)=d_H(y,x)$;
\item[3)] $d_H(x,y)d_H(y,z)\ge d_H(x,z)$.
\end{itemize}

In fact these properties show that $d_H$ is a semidistance
(recall that $\textbf{1}=0$ and $\odot=+$). Indeed, $d_H(x,y)=\textbf{1}=0$
whenever $x$ is equivalent to $y$ modulo 
\[
x\sim y\Leftrightarrow \exists\lambda\ne\textbf{0}: x=\lambda y. 
\] 
This semidistance is induced by the \emph{range seminorm}
\begin{equation}
\label{e:sn}
||x||=\bigoplus_{i,j\in K} x_ix_j^{-1},
\end{equation}
introduced in \cite{CG79}, see also~\cite{CunBut04}.
However, by a slight abuse of language we will refer to $d_H$ as 
to a distance. 

Now we measure the distance between an arbitrary $y\in \eig(A)$ and
the supporting plane $X_{ij}=\{\,x\mid x_ix_j^{-1}=A_{ij}\,\}$, i.e., the
minimal distance between $y$ and $x\in X_{ij}$. From~\cite{CGQ04}, Th.~18
it follows that this minimum is attained at the maximal vector of $X_{ij}$
not greater than $y$. Denote this vector by $y^{ij}$. Its
coordinates are very easy to find:
\begin{equation}
\label{e:yx}
\begin{array}{l}
y_l^{ij}=\bigoplus\{\, x_l\mid x_l\le y_l\,\}=y_l\quad\text{for $l\ne i,j$};\\[1.5ex]
y_i^{ij}=\bigoplus\{\, x_i\mid x_i\le y_i, A_{ij}^{-1}x_i\le y_j\,\}=A_{ij}y_j;\\[1.5ex]
y_j^{ij}=\bigoplus\{\, x_j\mid x_j\le y_j, A_{ij}x_j\le y_i\,\}=y_j.
\end{array}
\end{equation}

The distance~(\ref{e:hdist}) between $y$ and $X_{ij}$ is then equal to
\begin{equation}
\label{e:dyxij}
d_H(y,X_{ij})=A_{ij}^{-1}y_j^{-1}y_i.
\end{equation}

However, what we need is the
distance between $y$ and $\Gamma(\eig(A))$, i.e, the minimal distance
between $y$ and $\Gamma_{ij}$. The following proposition
makes our life simpler.
\begin{proposition}
\label{p:lifesimpler}
The distance between $y\in \eig(A)$ and the boundary
$\Gamma(\eig(A))$ is equal to the minimal distance between $y$ and
supporting planes.
\end{proposition} 

\textbf{Proof.} Clearly the minimal distance between $y$ and supporting 
planes is not greater than the distance between $y$ and $\Gamma(\eig(A))$.
Suppose $X_{ij}$ is the supporting plane such that the distance between
$y$ and $X_{ij}$ is minimal. This distance is equal to the distance
between $y$ and $y^{ij}$. If $y^{ij}$ belongs to $\eig(A)$ and hence to
$\Gamma_{ij}$ then we are done. Suppose not; then the system of
equalities $\bigoplus_l A_{kl}y_l^{ij}=y_k^{ij}$ must be violated
for some $k\in N$. Note that, if $k\ne i$, then $y_k^{ij}=y_k$
(see~(\ref{e:yx})), and there is no violation. So the violation must take
place for $k=i$. There must be an $l$ such that $A_{il}y_l^{ij}>y_i^{ij}$,
i.e. such that $A_{il}y_l>A_{ij}y_j$. Now consider $z$
such that $z_i=A_{il}y_l$ and $z_k=y_k$ for any $k\ne i$. Then $z$
belongs to $X_{il}$ and
$d_H(y,z)=A_{il}^{-1}y_l^{-1}y_i$. This distance is strictly less
than the distance between $y$ and $y^{ij}$, a contradiction.\hfill$\square$

Consequently,
\begin{equation}
\label{e:dygamma}
d_H(y,\Gamma(\eig(A)))=
\bigwedge_{i\ne j,A_{ij}\ne\textbf 0} A_{ij}^{-1} y_j^{-1}y_i.
\end{equation}

The key idea of Prop.~\ref{p:<>} below is that 
$\lambda(\A)^{-1}$, if $\lambda(\A)$ is invertible, is the largest
radius of Hilbert balls contained in $\eig(A)$. It can be said that 
$\lambda(\A)^{-1}$ is the radius of~\emph{inscribed} Hilbert balls,
as depicted on Fig.~5.

Let $\tau$ be any critical cyclic permutation of $\A$.
% The following
%two propositions provide information about the distances between
%$\eig(A_{\mu})$ and $\Gamma(\eig(A))$.
\begin{proposition}
\label{p:<>}
\begin{itemize}
\item[1)] In the case 
$\lambda(\A)>\textbf 0$ for any $y\in \eig(A)$ the distance
between $y$ and $\Gamma(\eig(A))$ is not greater than $\lambda(\A)^{-1}$;
%$d_H(y,\Gamma(\eig(A)))\le\lambda^{-1}$;
\item[2)] Let $\mu$ be such that $\lambda(\A)\le\mu<\textbf 1$, if 
$\lambda(\A)>\textbf 0$, or such that $\textbf{0}<\mu<\textbf 1$, if
$\lambda(\A)=\textbf 0$. Then for any
$y\in \eig(A_{\mu})$ the distance between $y$
and $\Gamma(\eig(A))$ is not less than $\mu^{-1}$;
\item[3)] In the case $\lambda(\A)>\textbf 0$, for any 
$i,j\in K(\tau)$ such that $j=\tau(i)$, 
and any $y\in \eig(A_{\lambda})$, the distance between $y$ and
the face $\Gamma_{ij}$ is equal to $\lambda(\A)^{-1}$.
%one has $d_H(y,\Gamma_{ij})=\lambda(\A)^{-1}$.
\end{itemize}
\end{proposition} 

\textbf{Proof.} 1) The distance between $y$ and $\Gamma(\eig(A))$ does
not exceed the minimal distance between $y$ and supporting planes $X_{ij}$
that correspond to the edges $(i,j)$ of the cyclic path determined by $\tau$, and
this minimal distance is not greater than $\lambda(\A)^{-1}$:
\[
\begin{array}{l@{{}\le{}}l}
d_H(y,\Gamma(\eig(A))) & \bigwedge_{i\in K(\tau)} A^{-1}_{i\tau(i)}y_{\tau(i)}^{-1}y_i\le\\[2ex]
&(\bigodot_{i\in K(\tau)} A_{i\tau(i)}^{-1})^{\frac{1}{\mid K(\tau)\mid}}=\lambda(\A)^{-1}.
\end{array}
\]

2)
If $y\in \eig(A_{\mu})$ then, since $\eig(A_{\mu})=\spann(A_{\mu}^*)$,
we have 
\begin{equation}
\label{e:y}
y=\bigoplus_{k\in M} \alpha_k (A_{\mu}^*)_{\tchk k},
\end{equation}
where $M=\supp(\alpha)$. Substituting~(\ref{e:y}) into~(\ref{e:dygamma}), we get
\begin{equation}
\label{e:dy}
d_H(y,\Gamma(\eig(A)))=
\bigwedge_{i\ne j,A_{ij}\ne\textbf 0} A_{ij}^{-1}
\bigwedge_{k\in M_j}\alpha_k^{-1} (A_{\mu}^*)_{jk}^{-1}
\bigoplus_l \alpha_l (A_{\mu}^*)_{il}.
\end{equation}
Here by $M_j$ we denote the set $M\meet \supp((A_{\mu}^*)_{j\tchk})$.
Now we estimate~(\ref{e:dy}) from below and use the inequalities
$A_{ij}\le\mu(A_{\mu}^*)_{ij}$ and $(A_{\mu}^*)_{ij} (A_{\mu}^*)_{jk}\le
(A_{\mu}^*)_{ik}$ (see (\ref{e:multid})):
\[
\begin{array}{l@{{}\ge{}}l}
d_H(y,\Gamma(\eig(A))) &
\bigwedge_{i\ne j,A_{ij}\ne\textbf 0}
\bigwedge_{k\in M_j} A_{ij}^{-1} (A_{\mu}^*)^{-1}_{jk} (A_{\mu}^*)_{ik}\ge\\[2ex]
& \bigwedge_{i\ne j,A_{ij}\ne\textbf 0}
\bigwedge_{k\in M_j}\mu^{-1} (A_{\mu}^*)_{ij}^{-1} (A_{\mu}^*)^{-1}_{jk} 
(A_{\mu}^*)_{ik}\ge\mu^{-1}.
\end{array}
\]

3) 
The distance between $y\in \eig(A_{\lambda})=\spann(A_{\lambda}^*)$ and
the supporting plane $X_{ij}$ is equal to
\begin{equation}
\label{e:dy1}
d_H(y,X_{ij})=
 A_{ij}^{-1}
\bigwedge_{k\in M_j}\alpha_k^{-1} (A_{\lambda}^*)_{jk}^{-1}
\bigoplus_{l\in M_i}\alpha_l (A_{\lambda}^*)_{il}.
\end{equation}

The cyclic permutation $\tau$ 
of $A_{\lambda}$ has the weight $\textbf 1$.
Hence for all $i,j\in K(\tau)$ such that $j=\tau(i)$ we have,
according to Prop.~\ref{p:proportia}, that 
$A_{ij}=\lambda(\A)(A_{\lambda}^*)_{ij}$ and
$(A_{\lambda}^*)_{ij}(A_{\lambda}^*)_{jl}=(A_{\lambda}^*)_{il}$.
Note that $(A_{\lambda}^*)_{ij}\ne\textbf 0$ for all $i,j\in K$
and therefore $(A_{\lambda}^*)_{jl}=\textbf 0$ if and only if
$(A_{\lambda}^*)_{il}=\textbf 0$, i.e. $M_i$ and $M_j$ coincide.
Making use of all this we write the upper estimate for $d_H(y,X_{ij})$:
\[
\begin{array}{l}
d_H(y,X_{ij}) \le
\bigoplus_{l\in M_i} A_{ij}^{-1}(A_{\lambda}^*)^{-1}_{jl} (A_{\lambda}^*)_{il}=\\[2ex]
=\bigoplus_{l\in M_i}\lambda(\A)^{-1} (A_{\lambda}^*)_{ij}^{-1} (A_{\lambda}^*)^{-1}_{jl} 
(A_{\lambda}^*)_{il}=\lambda(\A)^{-1}.
\end{array}
\]

We also have $d_H(y,\Gamma(\eig(A)))\ge\lambda(\A)^{-1}$ and therefore
(see Prop.~\ref{p:lifesimpler})
$d_H(y,X_{ij})=d_H(y,\Gamma_{ij})=\lambda(\A)^{-1}$.\hfill$\square$

The sets $\Gamma(\eig(A_{\mu}))$, for $\mu<\textbf 1$, 
are the subsets of $\eig(A)$
equidistant from
$\Gamma(\eig(A))$, as Prop.~\ref{p:=} suggests.

\begin{proposition}
\label{p:=}
For all $\mu$ such that $\lambda(\A)\le\mu<\textbf 1$, if 
$\lambda(\A)>\textbf 0$, or such that $\textbf 0<\mu<\textbf 1$,
if $\lambda(\A)=\textbf 0$, the distance 
$d_H(y,\Gamma(\eig(A)))$ is equal to $\mu^{-1}$ if and only if
$y\in\Gamma(\eig(A_{\mu})))$.
\end{proposition}

\textbf{Proof.} If $\lambda(\A)>\textbf 0$ and $\mu=\lambda(\A)$ 
then the statement readily follows
from the observation that $A_{\lambda}$ does not have strong permanent
and therefore (see Prop.~\ref{p:but}) $\eig(A_{\lambda})$ does not have
interior.

Let us consider $\mu>\lambda(\A)$. First, the equality
\[
\bigoplus_{i\ne j, A_{ij}\ne\textbf 0} A_{ij} y_jy_i^{-1}=\mu
\]
implies $A_{\mu}y=y$. So, if $d_H(y,\Gamma(\eig(A)))=\mu^{-1}$ then
$y\in \eig(A_{\mu})$. Assume that $y$ belongs to the interior of $\eig(A_{\mu})$.
Since $A_{\mu}$ is definite and has strong permanent, we can use
representation~(\ref{e:repbut}) or~(\ref{e:repbut1}) 
and obtain $\kappa<\textbf 1$
such that $y\in \eig(A_{\mu\kappa})$. Now statement 2) of 
Prop.~\ref{p:<>} implies that $d_H(y,\Gamma(\eig(A)))\ge(\mu\kappa)^{-1}>\mu$.
This is a contradiction, so $y\in\Gamma(\eig(A_{\mu}))$.

Suppose now that $y\in\Gamma(\eig(A_{\mu}))$. It means that there are $i\ne j$
such that $y_iy_j^{-1}=(A_{\mu})_{ij}$, where $(A_{\mu}^*)_{ij}=(A_{\mu})_{ij}$.
According to Prop.~\ref{p:injection}, this face corresponds to the
face of $\eig(A)$ determined by the entry $A_{ij}$, and the distance
between these two faces is clearly $\mu^{-1}$.\hfill$\square$

Throughout this section we dealt with eigenvectors having full
support. But let us recall the third statement of Prop.~\ref{p:defmat}. 
It says that
there might be eigenvectors with nontrivial support $K$.
The distance between these eigenvectors and part of any face with full support
would be infinite.  
Also, these eigenvectors
are eigenvectors of the submatrix $A_{KK}$. 
Therefore it is presumable,
in this case, to pose the problem of finding $d_H(y,\Gamma(\eig(A_{KK})))$.

We conclude this section with two examples.
\vskip1cm
\emph{Example 1} In the beginning of this section
we considered the definite matrix
\[
A=
\begin{pmatrix}
0&-4&1\\
1&0&1\\
-5&-7&0
\end{pmatrix},
\]
with the maximal cycle mean of $\A$ equal to $\lambda=-1.5$. Now
we pick the following three members of the $\{A_{\mu}\}$ family:
\[
A_{-0.5}=
\begin{pmatrix}
0&-3.5&1.5\\
1.5&0&1.5\\
-4.5&-6.5&0
\end{pmatrix},\;
A_{-1}=
\begin{pmatrix}
0&-3&2\\
2&0&2\\
-4&-6&0
\end{pmatrix},
\]

and
\[
A_{\lambda}=
\begin{pmatrix}
0&-2.5&2.5\\
2.5&0&2.5\\
-3.5&-5.5&0
\end{pmatrix}.
\]

%\begin{figure}[h]  
%\centering
%\vskip-5cm
%\includegraphics[width=13cm]{innstr}  
%\vskip-5cm
%\caption{The space $\spann(A)$ and the eigenspaces $\eig(A_{\mu})$}
%\end{figure}

\begin{figure}[htbp]
\label{f:amu}
\vskip0.5cm
\hskip1cm
\epsfysize=6cm
\epsffile{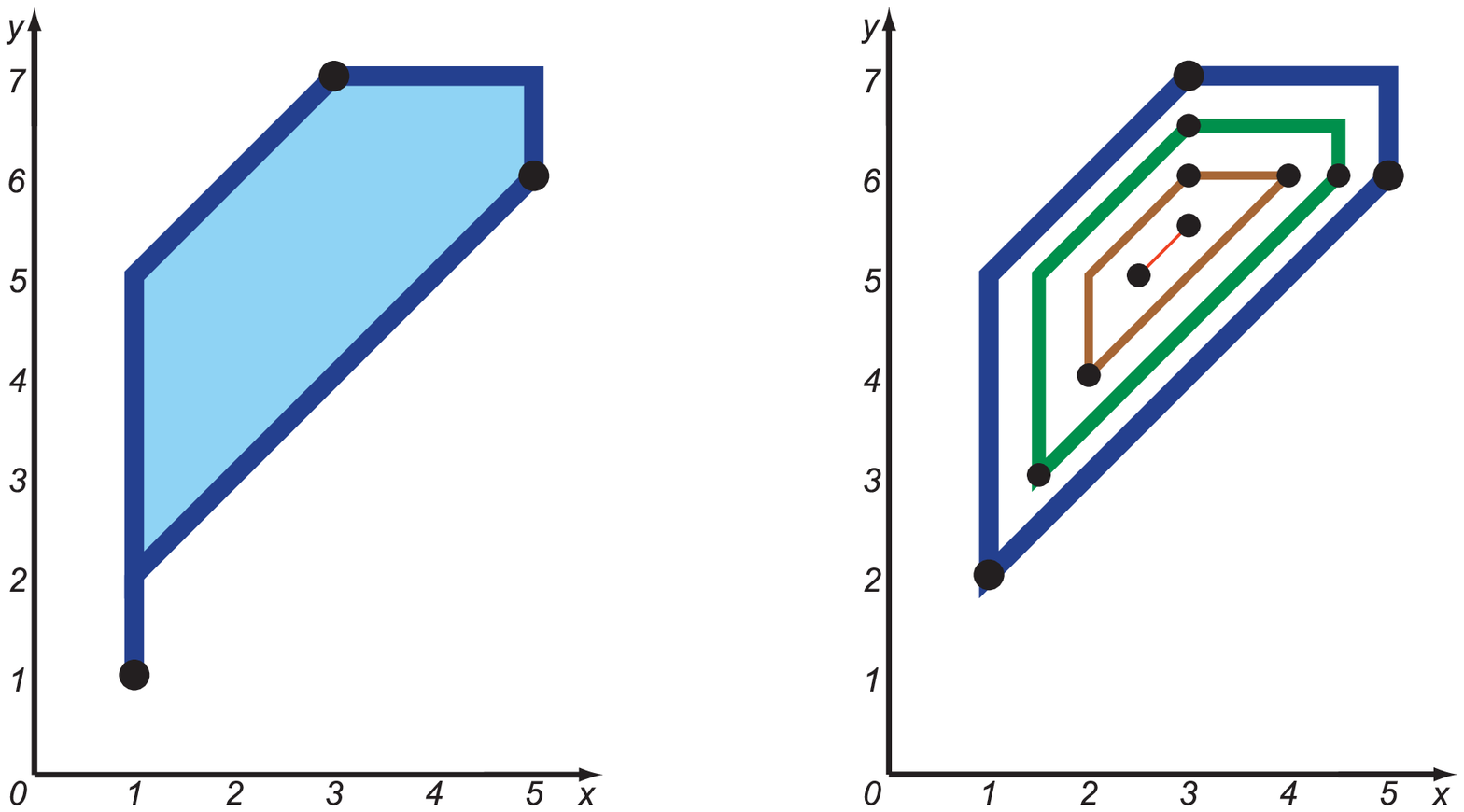}
\vskip0.5cm
\caption{The space $\spann(A)$ and the eigenspaces $\eig(A_{\mu})$}
\end{figure}

The left-hand side of Fig.~4 displays $\spann(A)$, and
the right-hand side displays the sets $\Gamma(\eig(A_{\mu}))$,
for $\mu=0$ (dark blue) $\mu=0.5$ (green), $\mu=1$ (brown), and 
$\mu=\lambda=1.5$ (red). The lines corresponding to
larger values of $\mu$ are given smaller weight (to help distinguishing
between different values of $\mu$ in black and white printing).
We see that there is an injection of systems
of inequalities describing $\eig(A_{\mu})$ into the system of
inequalities describing $\eig(A)$ in accordance with Prop.~\ref{p:injection}.

%\begin{figure}[h]  
%\centering
%\vskip-5cm
%\includegraphics[width=13cm]{hbinscr}  
%\vskip-5cm
%\caption{Two Hilbert balls inscribed in $\Gamma(\eig(A))$}
%\end{figure}

\begin{figure}[htbp]
\label{f:hballs}
\vskip0.5cm
\hskip1cm
\epsfysize=6cm
\epsffile{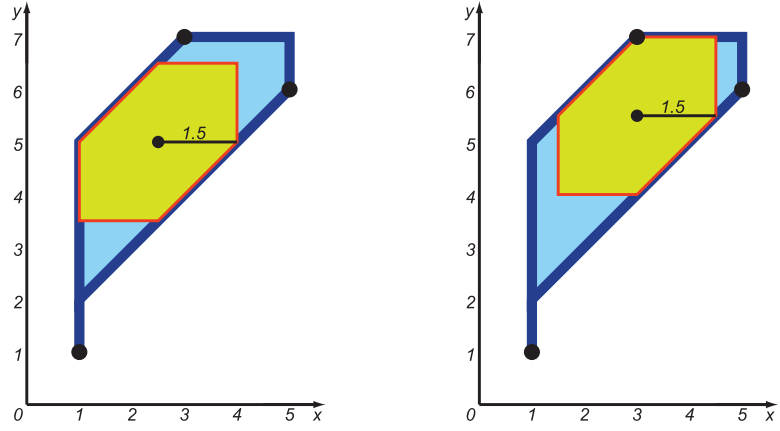}
\vskip0.5cm
\caption{Two Hilbert balls inscribed in $\Gamma(\eig(A))$}
\end{figure}

Fig.~5 displays, together with $\spann(A)$, two Hilbert
balls with the maximal radius $\lambda^{-1}=1.5$
inscribed in $\Gamma(\eig(A))$. The balls touch the `cycle faces'
$\{\,(x,y,z)\in \eig(A)\mid yx^{-1}=1\,\}$ and 
$\{\,(x,y,z)\in \eig(A)\mid xy^{-1}=-4\,\}$, 
in accordance with the third
statement of Prop.~\ref{p:<>}.
\vskip 1cm
\emph{Example 2} Consider a Hilbert ball with radius $d$ centered at $\{\lambda x\}$
($\lambda$ is any nonzero scalar).
It is the set
\[
Y=\{\;y\mid \bigoplus_{i,j} x_iy_i^{-1}y_jx_j^{-1}\le d\;\},
\]
or, equivalently, 
\[
Y=\{\;y\mid y_iy_j^{-1}\ge d^{-1}x_ix_j^{-1}\;\}.
\]
Denote by $D$ the matrix with entries $d_{ij}=d^{-1}x_ix_j^{-1}$. It is
easily verified that $D=D^*$. Then it follows from Prop.~\ref{p:syseq} that
the Hilbert ball is the eigenspace of $D$ and the columns of this matrix
are its generators. The maximal cycle mean of $\Tilde{D}$ is clearly
$d^{-1}$. The eigenspaces $\eig(D_{\mu})$ where $d^{-1}<\mu\le\textbf 1$
are Hilbert balls with radii $(\mu d)^{-1}$ centered at $\{\lambda x\}$,
and $\eig(D_{d^{-1}})$ is precisely $\{\lambda x\}$.

For a three-dimensional example, set $x=[5\;4\;0]^{T}$ and $d=3$. Then
\[
D=
\begin{pmatrix}
0&-2&2\\
-4&0&1\\
-8&-7&0
\end{pmatrix}.
\]
\newpage
%\begin{figure}[h]  
%\centering
%\vskip-4cm
%\includegraphics[width=7cm]{hballs}  
%\vskip-1cm
%\caption{Hilbert balls as eigenspaces}
%\end{figure}
\begin{figure}[htbp]
\label{f:hballs2}
\vskip-1cm
\hskip1cm
\epsfysize=7cm
\epsffile{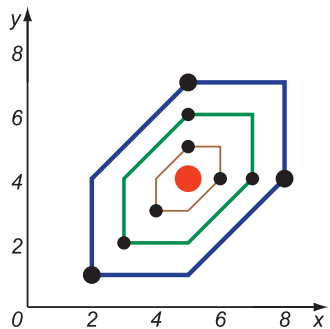}
%\vskip0.5cm
\caption{Hilbert balls as eigenspaces}
\end{figure}

Fig.~6 displays the sets $\Gamma(\eig(D_{\mu})$
for $\mu=0$ (dark blue), $\mu=-1$ (green), and $\mu=-2$ (brown) with
the same convention about the weight of lines as in Fig.~4. These
sets are concentric spheres centered at $\eig(D_{-3})=\{\lambda x\}$ 
(the large red circle in the center of Fig.~6).


\begin{thebibliography}{99}

\bibitem{BCOQ}
Baccelli~F., G.~Cohen, G.-J.~Olsder, and J.-P.~Quadrat,
\emph{Synchronization and Linearity. An Algebra for Discrete Event Systems.}
NY etc., Wiley, 1992.

\bibitem{But00}
Butkovi\v{c}~P.,
\emph{Simple Image Set of (max,+) Linear Mappings,}
Discrete Appl. Math. \textbf{105} (2000) 73-86.

\bibitem{But03}
Butkovi\v{c}~P.,
\emph{Max-algebra: Linear Algebra of Combinatorics?}
Linear Algebra and Appl. \textbf{367} (2003) 313-335.

\bibitem{Carre71}
Carr\'{e}~B.A.,
\emph{An Algebra for Network Routing Problems,}
J.~Inst.~Maths~Applics \textbf{7} (1971) 273-299.

%\bibitem{Frs}
%Cohen G., S.~Gaubert, and J.P.~Quadrat,
%Hahn-Banach Separation Theorem for Max-Plus Semimodules, in
%\emph{Optimal Control and Partial Differential Equations.}
%J.~L.~Menaldi, E.~Rofman and A.~Sulem Eds., IOS Press 2001.

\bibitem{CGQ04}
Cohen~G., S.~Gaubert, and J.-P.~Quadrat,
\emph{Duality and Separation Theorems in Idempotent Semimodules,}
Linear Algebra and Appl. \textbf{379} (2004) 395-422.
Also arXiv:math.FA/0212294.

\bibitem{CG79}
Cuninghame-Green~R.A.,
\emph{Minimax Algebra},
Lecture Notes in Economics and Mathematical Systems,
Vol. 166, Springer, 1979.

\bibitem{CunBut04}
Cuninghame-Green~R.A., and P.~Butkovi\v{c},
\emph{Bases in max-algebra},
Linear Algebra and Appl.\textbf{389} (2004) 107-120.

\bibitem{DSS}
Develin~M., F.~Santos, and B.~Sturmfels,
On the Rank of a Tropical Matrix. 
To appear in \emph{Discrete and Computational Geometry} 
(eds. J.E.~Goodman and J.~Pach),
MSRI Publications, Cambridge Univ. Press, 2005. Also arXiv:math.CO/0312114.

\bibitem{DevSt}
Develin~M., and B.~Sturmfels, 
\emph{Tropical Convexity},
Documenta Math. \textbf{9} (2004) 1-27. Also arXiv:math.MG/0308254.

\bibitem{Gau92}
Gaubert~S.,
\emph{Th\'{e}orie des Syst\`{e}mes Lin\'{e}aires dans les Dio\"{i}des},
Th\`{e}se, Ecole des Mines des Paris, Paris, 1992.

%\bibitem{Golan}
%Golan~J.,
%\emph{Semirings and Their Applications.}
%Kluwer, Dordrecht, 2000.

\bibitem{KolMas97}
Kolokoltsov~V.N, and V.P.~Maslov,
\emph{Idempotent Analysis and Its Applications,}
Kluwer Academic Publishers, Dordrecht, 1997.

\bibitem{LitMas95}
Litvinov~G.L. and V.P.~Maslov,
\emph{Correspondence principle for idempotent calculus
and some computer applications,}
Bures-Sur-Yvette: Institut
des Hautes Etudes Scientifiques
(IHES/M/95/33), 1995. See also: J.~Gunawardena(Editor),
\emph{Idempotency}, Publ. of the I.~Newton Institute,
Cambridge University Press, 1998, pp. 420-443.

\bibitem{LMS01}
Litvinov~G.L., V.P.~Maslov, and G.B.~Shpiz,
\emph{Idempotent Functional Analysis: an Algebraical Approach,}
Math. Notes \textbf{69}:5 (2001) 696-729.
Also arXiv:math.FA/0009128.

\bibitem{LitMaslova00}
Litvinov~G.L., and E.V.~Maslova,
\emph{Universal Numerical Algorithms and Their Software Implementation,}
Programming and Computer Software \textbf{26}:5 (2000) 275-280.


\bibitem{Mol88}
Moller~P.,
\emph{Th\'{e}orie Alg\'ebraique des Syst\`{e}mes \`{a} Ev\'{e}nements Discrets},
Th\`{e}se, Ecole des Mines des Paris, Paris, 1988.


\bibitem{Rote85}
Rote~G.,
\emph{A Systolic Array Algorithm for the Algebraic Path Problem,}
Computing \textbf{34} (1985) 191-219.

\bibitem{Wag88}
Wagneur~E.,
Moduloids and Pseudomodules.1.dimension theory. In \emph{Analysis
and Optimization of Systems} (eds. J.L.~Lions
and A.~Bensoussan), 
Lecture Notes in Control and Information Sciences, 1988.

\bibitem{Wag91}
Wagneur~E.,
\emph{Moduloids and Pseudomodules.1.dimension theory.}
Discrete Math. \textbf{98} (1991) 57-73.

\bibitem{Zimm81}
Zimmermann~U.,
\emph{Linear and Combinatorial Optimization in Ordered Algebraic Structures.}
North Holland, 1981.

\end{thebibliography}
\end{document}